\documentclass{article}

\usepackage[dvips]{graphicx}

\usepackage{amssymb,amsfonts,amsmath}

\usepackage{amsmath,amsmath,amsfonts}

\usepackage{amssymb}

\usepackage{epsfig}
\begin{document}

\def\s{\sigma}

\title{Implicit particle filters for data assimilation}

\author{Alexandre J. Chorin\\
Department of Mathematics\\
University of California at Berkeley\\ and Lawrence Berkeley National Laboratory\\
\\
Matthias Morzfeld\\
Department of Mechanical Engineering\\
University of California at Berkeley\\ and Lawrence Berkeley National
Laboratory\\
\\
Xuemin Tu\\
Department of Mathematics\\
University of California at Berkeley\\ and Lawrence
Berkeley National Laboratory}

\maketitle
\begin{abstract}

Implicit particle filters for data assimilation update the particles
by first choosing probabilities and then looking for particle locations 
that assume them, guiding the particles one by one to the high probability
domain. 
We provide a detailed description of these filters, with illustrative
examples, together with
new, more general, methods for solving the algebraic equations
and with a
new algorithm for parameter
identification. 
\end{abstract}

\section{Introduction}

{T}here are many problems in science, for example in meteorology and economics, in which the state
of a system must be identified from an uncertain equation supplemented by noisy data (see e.g. \cite{Do2,Stu}). A natural model of this situation consists of an Ito stochastic differential equation (SDE):
\begin{equation}
dx=f(x,t) \, dt+  g(x,t) \, dw,
\label{dynamics} 
\end{equation}
where $x=(x_1,x_2,\dots,x_m)$ is an $m$-dimensional vector, $f$ is an $m$-dimensional vector function, 
$ g(x,t)$ is an $m$ by $m$ matrix, and $w$ is Brownian motion which encapsulates all the uncertainty in the model. In the present paper we assume for simplicity that the matrix $g(x,t)$ is diagonal.  
The initial state
$x(0)$ is given and may be random as well.

The SDE is supplemented by measurements 
$b^n$ at times $t^n, n=0,1,\dots$.
The measurements are related to the state $x(t)$ by
\begin{equation}
b^n=h(x^n)+GW^n,
\label{observations}
\end{equation}
where $h$ is a $k$-dimensional, generally nonlinear, vector function with $k \le m$, $G$ is a matrix, $x^n=x(t^n)$,  and $W^n$ is a vector whose components are independent Gaussian variables of mean 0 and variance 1, independent also of the Brownian motion in equation (\ref{dynamics}). The
independence requirements can be greatly relaxed but will be observed in the present paper. The task  of a filter is to assimilate the data, i.e., estimate $x$ on the basis of both equation (\ref{dynamics}) and the observations (\ref{observations}).

If the system~(\ref{dynamics}) and equation (\ref{observations}) are linear and the data are Gaussian, the solution can be found in principle via the Kalman-Bucy filter (see e.g. \cite{Mil2}). In the general case, one often estimates $x$ as a statistic (often the mean) of a probability density function (pdf) evolving
under the combined effect of equations (\ref{dynamics}) and (\ref{observations}).
The initial state $x^0$ being known, all one has to do is
evaluate sequentially the pdfs $P_{n+1}$ of the variables $x^{n+1}$ given the equations and the data. 
In a ``particle" filter this is done by following 
``particles" (replicas of the system) whose empirical distribution at time $t^n$ approximates $P_n$. 
One may for example (see e.g.   
\cite{Ar2, Do2,Dou1,Mil2})
use the pdf $P_n$ and equation (\ref{dynamics}) to generate
a prior density (in the sense of Bayes) , and then use the data $b^{n+1}$ to generate sampling weights which define a
posterior density $P_{n+1}$. In addition, one has to sample backward to take into account the information each measurement provides about the past. This process can be very expensive because in most weighting schemes, most of the weights tend
to zero fast and 
the number of particles needed can grow catastrophically
(see e.g. 
\cite{Sny,Blb}); various strategies have been proposed to ameliorate this
problem. 

Our remedy is 
implicit sampling \cite{ct1,ct2}.
The number of particles needed in a filter remains moderate if one can find 
high probability particles;
to this end, 
implicit sampling works by first picking probabilities 
and then looking for particles that assume them, so that particles are guided efficiently
to the high probability region one at a time,  
without needing a global guess of the target density. In the present paper
we provide an expository account of particle filters,
separating clearly the general principles from details of 
implementation; we provide general solution algorithms for the resulting
algebraic equations, in particular for nonconvex cases which we had
not considered in our previous publications, as well as a new algorithm for
parameter identification based on an implicit filter. We also provide 
examples, in particular of nonconvex problems.  

Implicit filters are a special case of chainless/Markov field sampling
methods \cite{Ch101,Ch102}; a key connection was made in \cite{We1,We2}, where it was 
observed that in the sampling of stochastic differential equations, the marginals
needed in Markov field sampling can be read off the equations and need
not be estimated numerically.

\section{The mathematical framework}

The 
conditional probability density $P_n(x)$ at time $t^n$, determined by the SDE (\ref{dynamics})
given the observations (\ref{observations}), satisfies the recurrence
relation (see e.g. \cite{Do2}):
\begin{equation}
P_{n+1}(x^{n+1})=P_n(x^n)P(x^{n+1}|x^n)P(b^{n+1}|x^{n+1})/Z_0,
\label{general}
\end{equation}
where $P_{n+1}(x^{n+1})$ is the probability of the sample $x^{n+1}$ at time $t^{n+1}$ given the observations $b^j$ for $j\le {n+1}$,
$P_n(x^n)$ is the probability of a the sample $x^n$ at time $t^n$ given the
observations $b^j$ for $j\le n$,
$P({x}^{n+1}|x^n)$ is the probability of a sample $x^{n+1}$ at time $t^{n+1}$
given a sample $x^n$ at time $t^n$,  
$P(b^{n+1}|x^{n+1})$ is the probability of the observations $b^{n+1}$ given
the sample $x^{n+1}$ at time $t^{n+1}$, and $Z_0$ is a normalization constant
independent of $x^n$ and $x^{n+1}$. This is Bayes' theorem. 

We estimate $P_{n+1}$ with the help of M particles, with positions $X_i^n$
at time $t^n$ and $X^{n+1}_i$ at time $t^{n+1}$ ($i=1,\dots,M$), which define
empirical densities $\hat{P}_n,\hat{P}_{n+1}$ that approximate  
$P_n,P_{n+1}$. We do this by requiring that, when a particle moves from $X_i^n$
to $X^{n+1}_i$ the probability of $X^{n+1}_i$ be 
\begin{equation}
 P(X^{n+1}_i)=  P(X^n_i)P(X^{n+1}_i|X^n_i)P(b^{n+1}|X^{n+1}_i)/Z_0, \label{bayes}
\end{equation}
where the hats have been omitted, 
$P(X^n_i)$, the probability of $X^n_i$, is assumed given, 
the pdf $P(X^{n+1}_i|X^n_i)$, the probability of $X^{n+1}_i$ given $X^n_i$, is determined by the SDE (\ref{dynamics}), the pdf $P(b^{n+1}|X^{n+1}_i)$, the probability of the observations $b^{n+1}$ given the new positions $X^{n+1}_i$, is determined by the observation equation (\ref{observations}). 
We shall see below that one can set $P(X^n_i)=1$ without loss of generality.

Equation (\ref{bayes}) defines the pdf we need to sample for each particle; this pdf is
known, in the sense that once one has a sample, one can evaluate its
probability (up to a constant); the difficulty is to find high probability samples,
especially when the number of variables is large. 
The idea in implicit sampling is to define probabilities first, and then
look for particles that assume them; this is done by  
choosing once and for all a fixed reference random variable, say $\xi$, with a given
pdf, say a Gaussian $\exp(-\xi^T\xi/2)/(2\pi)^{m/2})$, which one knows 
how to sample so that most samples have high probability, and then
making $X^{n+1}_i$ a function of $\xi$, a different function of
each particle and each step, each function designed so that the map $\xi \rightarrow
X^{n+1}_i$ connects highly probable values of $\xi$ to 
highly probable values of $X^{n+1}_i$. To that end, 
write $$P(X^{n+1}_i|X^n_i)P(b^{n+1}|X^{n+1}_i)=\exp(-F_i(X)),$$ where on the right-hand
side $X$ is a shorthand for $X^{n+1}_i$ and all the other arguments
are omitted. This defines a function $F_i$ for each particle $i$ and each time $t^n$. For each $i$ and $n$, $F_i$ is an explicitly known function of $X=X^{n+1}_i$. Then solve the
equation 
\begin{equation}
F_i(X)-\phi_i=\xi^T\xi/2,
\label{main}
\end{equation}
where
$\xi$ is a sample of the fixed reference variable and 
$\phi_i$ is an additive factor needed to make the equation solvable. 
The need for $\phi_i$ becomes obvious if one considers the case of a linear 
observation function $h$ in equation (\ref{observations}), so that the right side
of equation (\ref{main}) is quadratic but the left is 
a quadratic plus a constant. It is clear that setting $\phi=\min F$ will
do the job, but this is not necessarily the best choice (see below). 
We also require that for each particle, the function $X^{n+1}_i=X=X(\xi)$ defined by
(\ref{main}) be one-to-one
so that the correct pdf is sampled, in particular, it must have distinct branches for $\xi >0$
and $\xi <0$. 
The solution of 
(\ref{main}) is discussed in the
next section.  
From now on we omit the index $i$ in both $F$ and $\phi$, but it
should not be forgotten that these function vary from particle to particle
and from one time step to the next.

Once the function $X=X(\xi)$ is determined, each value of $X^{n+1}=X$ (the subscript $i$ is omitted)  appears with
 probability $\exp(-\xi^T\xi/2)J^{-1}/(\pi)^{m/2}$, where $J$ is the Jacobian of the
map $X=X(\xi)$, 
while the product $P(X^{n+1}|X^n)P(b^{n+1}|X^{n+1})$ evaluated
at $X^{n+1}$ equals $\exp(-\xi^T\xi/2)\exp(-\phi)/(2\pi)^{m/2}$. 
The sampling weight for the particle is therefore $\exp(-\phi)J$. If the map $\xi \rightarrow X$ is smooth near $\xi=0$, so that 
$\phi$ and $J$ do not vary rapidly from particle to particle, and if   
there is an easy way to compute $J$ (see the next section), then we have an effective way to
sample $P_{n+1}$ given $P_n$.
It is important to note 
that though the functions $F$ and $\phi$ vary from particle to particle, the probabilities
of the various samples are expressed in terms of the fixed reference pdf,
so that they can be compared with each other. 

The weights can be eliminated by resampling. 
A standard resampling algorithm goes as follows \cite{Do2}: let the weight of the
$i$-th particle be $W_i, i=1,\dots,M$. Define $A=\sum W_i$;
for each of $M$ random numbers
$\theta_k,k=1,\dots,M$ drawn from the uniform distribution on $[0,1]$,
choose a new ${\widehat X}^{n+1}_k=X^{n+1}_i$
such that
$A^{-1}\sum_{j=1}^{i-1}W_j <  \theta_k \le A^{-1}\sum_{j=1}^i W_j$,
and then suppress the hat. This justifies the statement following equation 
(\ref{bayes}) that one can set $P(X_n)=1$.

To see what has been gained, compare our construction with the 
usual ``Bayesian" particle filter, where one samples
$P(X^{n+1}|X^n)P(b^{n+1}|X^{n+1})$ by first finding a ``prior" density
$Q(X^{n+1})$ (omitting all arguments other than $X^{n+1}$), such that the ratio $W=P(X^{n+1})/Q(X^{n+1})$ is close to a constant, and then assigning to the $i$-th particle 
the importance weight $W=W_i$ evaluated at the location of the particle. 
The pdf defined by the
set of positions and weights is the density $P_{n+1}$ we are looking for.
An important special case 
is the choice $Q(X^{n+1})=P(X^{n+1}|X^n)$; the prior is then defined by the equation
of motion alone and the posterior is obtained by using the observations to
weight the particles. We shall refer to this special case as ``standard importance
sampling" or ``standard filter". 
Of course, once the positions and the weights of the particles have been determined,
one should resample as above. 

The catch in these earlier constructions is that the prior density $Q$ and the desired posterior 
can come close to being  
mutually singular, and the number of particles needed may become
catastrophically large, especially when the number of variables $m$ is large.
To avoid this catch one has to make a good guess for the pdf $Q$,
which may not be easy because $Q$ should approximate the density $P_{n+1}$
one is looking for- this is the basic conundrum of Monte Carlo methods,
in which one needs a good estimate to get a good estimate.
In contrast, in 
implicit sampling one does a separate calculation for each sample and there is
no need for prior global information. 
One can of course still identify
the pdf defined by the positions of the particles at time $t^{n+1}$ as a ``prior"
and the pdf defined by both the positions and the weights as
a ``posterior" density. 

Finally, implicit sampling can be viewed as an implicit Monte Carlo scheme for solving the
Zakai equation \cite{Zak}, which  describes the evolution of the (unnormalized) conditional distribution 
for a SDE conditioned by observations. This should be contrasted with
the procedure in the popular ensemble Kalman filter (see e.g. \cite{Even}), where a Gaussian approximation
of the pdf defined by the SDE is extracted from a Monte Carlo solution
of the corresponding Fokker-Planck equation, a Gaussian approximation is
made for the pdf $P(b^{n+1}|x^{n+1})$, and new particle positions
are obtained by a Kalman step. Our replacement of the Fokker-Planck 
equation that corresponds to the SDE alone by a Zakai equation that
describes the evolution of the unnormalized conditional distribution
does away with the need for the approximate and expensive extraction
of Gaussians and consequent Kalman step.

\section{Solution of the algebraic equation that defines a new sample}

We now explain how to solve equation (\ref{main}),
$F(X)-\phi=\xi^T\xi/2$, under 
several sets of assumptions which are met in practice. 
Note the great latitude this equation provides in linking the $\xi$
variables to the $X$ variables;
equation (\ref{main}) is a single equation
that connects $2m$ variables (the $m$ components of $\xi$ and the $m$ components
of $X$) and can be satisfied by many maps $\xi \rightarrow X$; these are useful as long as
(i) they are one-to-one, (ii) they map the neighborhood of $0$ into a set that contains the minimum of $F$, (iii) they are smooth near $\xi=0$ so that the weights $\exp(-\phi)$ and the Jacobian $J$ not vary unduly from particle to particle in the
target area, and (iv) they allow the Jacobian $J$ to be calculated easily. The solution methods presented here are far from exhaustive; further examples and refinements
in the context of specific applications.

{\bf
Algorithm (A)} (presented in \cite{ct1,ct2}) : Assume the function $F$ is convex upwards.
For each particle, we set up an iteration, with iterates $X^{n+1,j}$, $j=0,1,\dots$, ($X^j$ for brevity), with $X^0=0$, 
that converge to 
the next position $X^{n+1}$ of that particle.
The index $i$ that identifies the particle is omitted again.
We write the equations as if the system were one-dimensional;
the multidimensional case was presented in detail in \cite{ct2}. First we sample the reference
variable $\xi$. The iteration is defined when one knows how to 
find $X^{j+1}$ given $X^j$. 

Expand the observation function $h$ in equation 
(\ref{observations}) around $X^j$: 
\begin{equation}\label{lineeqn}
h(X^{j+1})=h(X^j)+(Dh)^j(X^{j+1}-X^j),
\end{equation}
where $(Dh)^j$ is the derivative of $h$ evaluated at $X^j$.
The observation equation (\ref{observations}) is now approximated as
a linear function of $X^{j+1}$, and the function $F$ is the sum of two
Gaussians in $X^{j+1}$. Completing a square yields a single Gaussian with a remainder
$\phi$, i.e., 
$F(X)=(x-{\bar a})^2/(2{\bar v})+\phi(X^j)$, 
where the parameters $\phi,{\bar a},{\bar v}$  are functions of $X^j$
(this is what we called in \cite{ct1}  a ``pseudo-Gaussian").
The next iterate is now $X^j=
{\bar a}+\sqrt{\bar v}\xi$. 
In the multidimensional case, when each component of the function $h$ in equation (\ref{observations}) depends on more than one variable, finding $X$ as a function of $\xi$ 
may require the solution of a linear system of equations, which can be performed e.g. by a Choleski factorization, as in \cite{ct2},
or by a rotation, as in \cite{ct1}. 
If the iteration converges, it converges to 
the exact solution of equation (\ref{main}), with $\phi$ the limit of the
$\phi(X^j)$. Its convergence can be accelerated by Aitken's extrapolation \cite{Isa}. The Jacobian $J$ can be evaluated either by an implicit
differentiation of equation (\ref{main}) or numerically, by  
perturbing $\xi$ in equation (\ref{main}) and solving the perturbed equation
(which should not require more than a single additional iteration step).
It is easy to see that this iteration, when it converges, produces a mapping $\xi \rightarrow X$
that is one to one and onto.

An important special case occurs when the observation function $h$ is linear in $X$;
it is immaterial  
whether the SDE ({\ref{dynamics}) is linear. In this case
the iteration converges in one step;
the Jacobian $J$ is easy to find; if in addition the function $g(x,t)$ in equation (\ref{dynamics}) is independent of $x$, then $J$ is independent of the particle 
and need not be evaluated;
the additive term $\phi$ can be written explicitly
as a function of the previous position $X^n$
of the particle and of the observation
$b^{n+1}$.
We recover an easy implementation of optimal 
sequential importance sampling (see e.g. \cite{Ar2,Do2,Dou1}).

This iteration has been used in \cite{ct2}. It may
fail to converge if the function $F$ is not convex (as 
happens in particular when the observation function $h$ is highly nonlinear).
One may resort then to
the next construction.

{\bf
Algorithm (B)}. Assume the function $F$ is $U$-shaped, i.e., in the scalar case, it is at least piecewise differentiable, $F'$ vanishes at a single point which is 
a minimum, $F$ is strictly decreasing on one side of the minimum and
strictly increasing on the other, with 
$F(X)=\infty$ when $X=\pm \infty$. In the $m$-dimensional case, assume that $F$
has a single minimum and that each intersection of the graph of the function $y=F(X)$
with a vertical plane through the minimum is $U$-shaped in the scalar sense 
(note that a function may be $U$-shaped without being convex).

Find $z$, the minimum of $F$ (note that this is the minimum of a given real valued
function, not a minimum of a possibly multimodal pdf generated by the SDE; finding
this minimum is not equivalent to the difficult problem of finding a maximum likelihood estimate of the state
of the system).  The minimum $z$ can be found by standard minimization algorithms.

Again we are solving the equations by finding iterates $X^j$ that converge to $X^{n+1}$. 
In the scalar case, given a sample of the reference variable $\xi$, find first $X^0$
such that $X^0-z$ has the sign of $\xi$, and then find the next iterates $X^j$ by standard tools (e.g. by Newton iteration), modified so
that the $X^j$ are prevented from leaping over $z$.  

In the vector case, if the observation function is diagonal, i.e. 
each component of the observation is a function of a single component of the solution $X$, then 
the scalar algorithm can be used component by component. In more complicated
situations one can take advantage of 
the freedom in connecting $\xi$ to $X$. 

Here is an interesting example of the use of this freedom, which we present in the case of a multidimensional
problem where the observation function is linear but need not be diagonal. Set $\phi=\min F$. The function $F(X)-\phi$ can 
now be written as $(X-a)^TA(X-a)/2$, where $a$ is a known vector, $T$ denotes a transpose as before, and $A$ is a positive
definite symmetric matrix. Write further $y=X-a$. Equation (\ref{main}) becomes
\begin{equation}
y^TAy=|\xi|^2,
\label{cur1}
\end{equation}
where $|\xi|$ is the length of the vector $\xi$. Make the ansatz:
$$y=\lambda\eta,$$
where $\lambda$ is a scalar, $\eta=\xi/|\xi|$ is a random unit vector and 
$\xi$ is a sample of of the reference density. 
Substitution into (\ref{cur1}) yields
\begin{equation}
\lambda^2(\eta^TA\eta)=|\xi|^2.
\label{cur2}
\end{equation}

It is easy to see that $E[\eta_i\eta_j]=\delta_{ij}/m$, where $E[\cdot]$ denotes
an expected value, the $\eta_i$ are the components of $\eta$, $m$ is the number of variables, and $\delta_{ij}$
is the Kronecker delta, and hence
$$E[\eta^TA\eta]={\rm trace }(A)/m.$$
Replace equation (\ref{cur2}) by 
\begin{equation}
\lambda^2\Lambda=|\xi|^2.
\label{cur3}
\end{equation}
where $\Lambda={\rm trace }(A)/m$.
This equation has the solution $\lambda=|\xi|/\sqrt{\Lambda}$,
and substitution into the ansatz leads to $y_i=\xi_i/\sqrt{\Lambda}$,
a transformation with Jacobian $J=\Lambda^{-m/2}$. 
The difference between equations (\ref{cur2}) and (\ref{cur3})
can be compensated for by adding to $\phi$ the term
$\lambda^2[(\eta^TA\eta)-\Lambda]$.
Notice now that as
$m\rightarrow \infty, (\eta^TA\eta)\rightarrow \Lambda$ (a stochastic weak
law of large numbers), so that when the number of variables is sufficiently large, the perturbation one has to compensate for becomes
negligible. Generalizations and applications of this construction will be given elsewhere in the context of specific applications.

One can readily devise algorithms also for cases where $F$ is not $U$-shaped, for example,
by dividing $F$ into monotonic pieces and sampling each of these pieces with its predetermined
probability. An alternative that is usually easier is to replace the non-$U$-shaped function $F$ by a suitable 
$U$-shaped function $F_0$ and make up for the bias by adding $F_0(X)-F(X)$  to
the additive term $\phi$; see the examples below.

\section{Backward sampling and sparse observations}

The algorithms of the previous sections are sufficient to create a  
filter, but accuracy may require an additional step. Every observation provides information
not only about the future but also about the past- it may, for example,
tag as improbable earlier states that had seemed probable
before the observation was made; in general one has to 
go back and correct the past after every observation (this backward sampling is often misleadingly motivated solely by the need to create greater diversity
among the particles in a Bayesian filter). A detailed construction has been
presented in \cite{ct2}; the examples in the present paper are simple
enough so that backward sampling does not significantly enhance their
performance, so we will be content here with presenting the construction
in principle, without much detail; it is a straightforward extension 
of the work above.

Consider the $i$-th particle, and suppose we have sampled 
its positions $X^{n-1}$, $X^n,X^{n+1}$ at times $t^{n-1},t^n,t^{n+1}$. 
Now we would like to go back and resample a new position $X^n$
at time $t^n$, given $X^{n-1}$ and $X^{n+1}$. The probability
density of $X=X^n$ is proportional to $P(X)=P(X|X^{n-1})P(b^n|X)P(X^{n+1}|X)$.
Write $P(X)=\exp(-F(X))$, sample a Gaussian reference variable $\xi$,
solve $F(X)-\phi=\xi^T\xi/2$ as above, and you are done. If need be,
one can then go further back and resample $X^{n-1}, X^{n-2},\dots$
Note that backward sampling relates $P_{n+1}$ to $P_{n-k}$ for
$k \ge 0$.

A similar construction can be used when the observations are sparse in time,
for example, if the time step needed to discretize the SDE accurately
is shorter than the time interval between observations.
Suppose we have sampled $X^{n-1}$, have an observation at time $t^{n+1}$
but not at time $t^n$, so that we have to sample simultaneously 
$X^n$ and $X^{n+1}$ from the SDE and the observation $b^{n+1}$. The
joint probability of $X=(X^n,X^{n+1})$ is proportional to
$P(X^n|X^{n-1})P(X^{n+1}|X^n)P(b^{n+1}|X^{n+1})$. Again, write this probability as $\exp(-F(X))$
and equate $F(X)-\phi$ to $\xi^T\xi/2$, where $\xi$ is a $2M$-dimensional reference variable.
Detailed expression for the vector case, as well as examples, can be found 
in \cite{ct2}.

\section{Examples}

We now present examples that illustrate the algorithms we
have just described. For more examples, see \cite{ct1,ct2}. 

We begin with a response to a comment we have often heard:
"this is nice,
but the construction will fail the moment you are faced with potentials with 
multiple wells". This is not so- the function $F$ depends on the nature of the noise in the SDE and on the function $h=h(x)$ in the observation equation (\ref{observations}), but not on the potential. Consider for example a one dimensional particle moving in
the potential 
$ V(x)=2.5(x^2-0.5)^2$, (see Figure 1), with the force 
$f(x)=-\nabla V=-10x(x^2-1)$ and the resulting SDE  
$dx=f(x)dt+\sigma dw$, where $\sigma=.1$ and $w$ is Brownian motion with unit variance;
with this choice of parameters the SDE has an invariant density concentrated
in the neighborhoods of $x=\pm \sqrt{1/2}$. 
We consider 
linear observations $b^n=x(t^n)+W$, where $W$ is a Gaussian
variable with mean zero and variance $s=.025$. We approximate the SDE by an Euler
scheme \cite{Klo1} 
with time step $\delta=0.01$, and assume observations are available
at all the points $n\delta$.
The particles all start at $x=0$. We produce data $b^n$
by running a single particle and adding to its positions errors drawn from the assumed error density in equation (\ref{observations}),
and then attempt to reconstruct this path with our filter.

\begin{figure}
\centerline{
{\includegraphics[width=.8\textwidth]{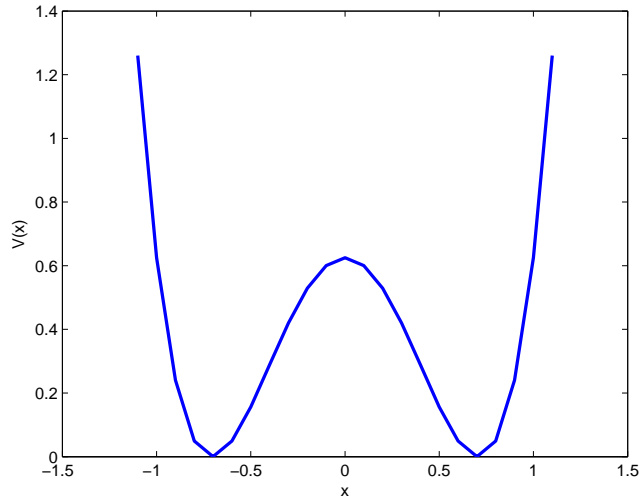}}}
\caption{The potential in the first example.}
\end{figure}

\begin{figure}
\centerline{
{\includegraphics[width=.8\textwidth]{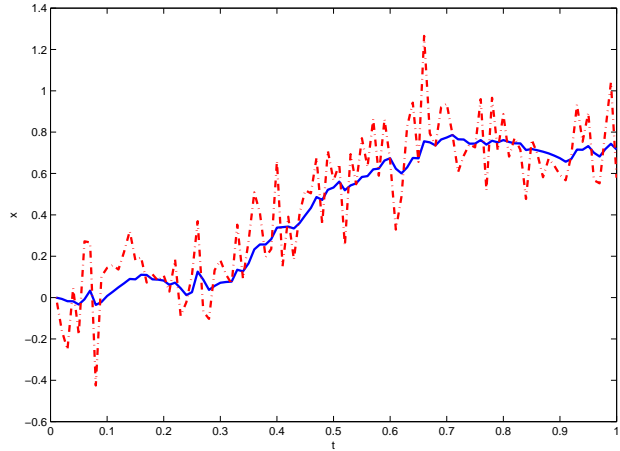}}}
\caption{A random path (broken line) and its reconstruction by our filter (solid line).}
\end{figure}

For the $i$-the particle located at time $n\delta$ at $X^n_i$ the function $F(X)$ is
$$F(X)=(X-X^n_i)^2/(2\sigma\delta)+(X-b^{n+1})^2)^2/(2s),$$ which is
always convex. A completion of a square yields $\min F=\phi=(1/2)(X^n_i-b^{n+1})^2/(\sigma\delta+s)$;
the Jacobian $J$ is independent of the particle and need not be evaluated. In Figure
2 we display a particle run used to generate data and its reconstruction by our filter
with $50$ particles. This figure is included for completeness but 
both of these paths are random, their difference varies from realization to
realization, and may be large or small by accident. To get a quantitative 
estimate of the performance of the filter, we repeated this calculation $10^4$ times
and computed the mean and the variance of the difference $\Delta$ between
the run that generated the data and its reconstruction at time $t=1$, see Table I. 
This Table shows that the filter is unbiased and that the variance of $\Delta$
is comparable to the variance of the error in the observations $s=0.025$.
Note that even with one single particle (and therefore no resampling)
the results
are still acceptable.

\begin{table}
\begin{center}
Table I
\end{center}
Mean and variance of the discrepancy between the observed path
and the reconstructed path in example 1 as a function of the
number of particles M, with $s=0.025$.
\begin{center}
\begin{tabular}{ccc}
M &     mean &    variance\\
100 &   -.0001  &  .021\\
50 &   -.0001 &   .022\\
20 &   -.0001 &   .023\\
10 &    .0001 &     .024\\
5  &   -.0001 &    .027\\
1  &   -.0001 &    .038
\end{tabular}
\end{center}
\end{table}

We now discuss the relation between the posterior we wish to sample and the
prior in several special cases, including non-convex situations. 
We want to produce samples of the pdf
$P(x)=\exp(-F(x))/Z$, where $Z$ is a normalization constant and
\begin{equation}
F(x)=x^2/(2\sigma)+(h(x)-b)^2/(2s)
\label{rare}
\end{equation}
and $h(x)$ is a given function of $x$ (as in equation (\ref{observations})) and $\sigma,s,b$ are given parameters. This can be viewed as a the first time 
step in time for a filtering problem where all the particles
start from the same point so that $\exp(-F(x))/Z=P_1$, or as an analysis
of the sampling for one particular particle in a general filtering problem,
or as an instance of the more general problem of sampling a given pdf
when the important events may be rare. In standard Bayesian
sampling one samples the variable with pdf $\exp(-x^2)/(2\sigma))/\sqrt{2\pi\sigma}$ and
then one attaches to the sample at $x$ the weight $\exp(-(h(x)-b)^2/(2s))$;
in an implicit sampler one finds a sample $x$ by solving $F(x)-\phi=\xi^2/2$
for a suitable $\phi$ and $\xi$ and attaching to the sample the weight
$\exp(-\phi)J$. For given $\sigma, s$, the problem becomes more
challenging as $|b|$ increases.

In both the standard and the implicit filters 
one can view the empirical pdf generated by the unweighted samples
as a ``prior" and the one generated by the weighted samples as the ``posterior". The difficulty with standard filters is that the prior and posterior densities may approach being mutually singular,  
so it is of interest to estimate the Radon-Nikodym derivative of one of these with respect to the other. If that derivative is a constant, we have achieved perfect importance
sampling, as every neighborhood in the sample space is visited with a frequency
proportional to its density.
We estimate the Radon-Nikodym derivative of the prior with respect to the posterior
as follows.
In this simple problem one can evaluate the probability of any interval
with respect to the posterior we wish to sample by quadratures.
We divide the interval $[0,1]$ into $K$ pieces of equal lengths $1/K$, then find numerically points $Y_1,Y_2,\dots, Y_{K-1},$ with $Y_K=+\infty$, such
that the posterior probability of the interval $[-\infty,Y_k]$ is
$k/K$ for $k=1,2,\dots,K$.
We then find $L=10^5$ samples of the prior and plot of a histogram of the frequencies
with which these samples fall into the posterior equal probability intervals
$(Y_{k-1},Y_k)$. The more this histogram departs from being a constant
independent of $k$, the more samples are needed to calculate the
statistics of the posterior.

If $h(x)$ is linear, the weights in the implicit filter are all equal
and the histogram is constant for all values of $b$. 
 This remains true for all values of $b$,
i.e., however far the observation $b$ is from what one may expect from the SDE alone.
This is not the case with a standard Bayesian filter, where some parts of the sample space
that have non-zero probability are visited very rarely.
In Table II we list the histogram of frequencies for a linear observation
function $h(x)=x$ and $b=2$ in a standard Bayesian filter, with K=10. We used $10^4$ samples; 
the fluctuations in the implicit case measure only the accuracy with which the
histogram is computed with this number of samples. 

\begin{table}
\begin{center}
Table II
\end{center}
Histogram of the Radon-Nikodym derivative of the prior with respect to
the posterior, standard Bayesian filter vs. the implicit filter,
10000 particles, $b=2$, $\sigma=s=0.1$, $h(x)=x$.
\begin{center}
\begin{tabular}{ccc}
k&        standard&   implicit \\
1&        .987&        .099\\
2 &       .006&        .108\\
3&        .002&       .097\\
4 &       .001&        .099\\
5 &       .004&        .101\\
6 &       .003&        .099\\ 
7 &       .001&        .101\\
8 &       .001 &       .101\\
9 &       .000 &       .102\\
10 &      .000  &      .093
\end{tabular}
\end{center}
\end{table}
As a consequence, estimates obtained with the implicit filter are much more reliable 
than the ones obtained with the standard Bayesian filter. In Table III we list the
estimates of the mean position of the linear problem as a function of b, with 30 particles,
$\sigma=s=0.1$,
for the standard Bayesian and the implicit filters, compared with the exact result.
The standard deviations are not displayed, they are all near 0.01.

\begin{table}
\begin{center}
Table III
\end{center}
Comparison of the the estimates of the means, implicit vs.\ standard filter,
$30$ particles, together with the exact results, linear case,
as explained in the text.
\begin{center}
\begin{tabular}{cccc}
b &     exact&     standard&  implicit\\
0&      0&          -.05&      .02\\
0.5 &   .25 &        .10  &    .27\\
1.&     .5&          .18   &   .51\\
1.5 &   .75&         .23&      .76\\
2.&     1.&          .26&      1.01
\end{tabular}
\end{center}
\end{table}
The results in this one-dimensional problem mirror the situation with the
example of Bickel et al. \cite{Blb,Sny}, designed to display the breakdown of the
standard Bayesian filter when the number of dimension is large; what happens there
is that one particle hogs almost the whole weight, so that the number of particles
needed grows catastrophically; in contrast, the implicit filter assigns equal weights
to all the particles in any number of dimensions, so that the number of particles
needed is independent of dimension, see also \cite{ct2}.

\begin{figure}
\centerline{
{\includegraphics[width=.8\textwidth]{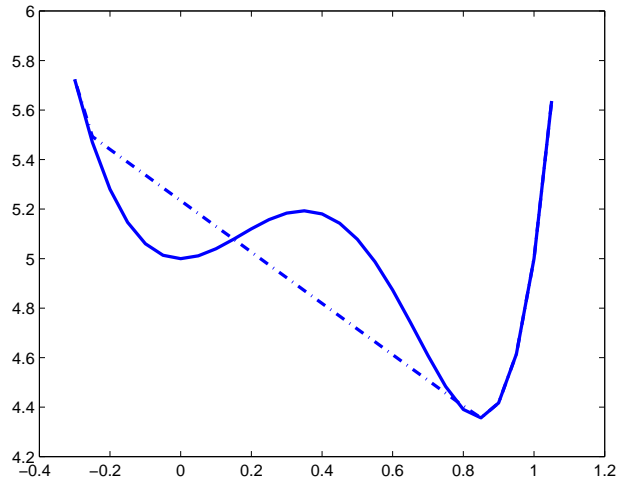}}}
\caption{A non-convex function $F$ (solid line) and a $U$-shaped substitute (broken line).}
\end{figure}

We now turn to nonlinear and nonconvex examples. Let the observation function
$h$ be strongly nonlinear: $h(x)=x^3$. With $\sigma=s=0.1$; the pdf (\ref{rare}) becomes
non-$U$-shaped for $|b| \ge .77$.  In Figure 3 we display the function $F$
for $b=1$ (the solid curve). To use the algorithms above we need a
substitute function $F_0$ that is $U$-shaped; we also display in Figure 3
(the broken line)
the function $F_0$ we used; the recipe here is to link a point above the 
local minimum on the left to the absolute minimum on the right by a straight
line. There are many other possible constructions; the only general rule
is to make the minimum of $F_0$ equal the absolute minimum of $F$, for obvious reasons.
As described above, we solve $F_0(x)-\phi=\xi^T \xi/2$ and set
$\phi=\min F_0+F_0(x)-F(x)$. 
It is important to note that this construction does not introduce any bias.
The function $F_0$ constructed in this way is $U$-shaped but need not be convex, so that one needs 
algorithm (B) described above. 
In Table IV  we compare the Radon-Nikodym derivatives of the prior 
with respect to the posterior for the resulting implicit sampling and for
standard Bayesian sampling with $\sigma=s=0.1, b=1.5$. 

\begin{table}
\begin{center}
Table IV
\end{center}
Radon-Nikodym derivatives of the prior with respect to the posterior,
$h(x)=x^3, \sigma=s=0.1, b=1.5,$ $10000$ samples,  $F_0$ as in the text.
\begin{center}
\begin{tabular}{ccc}
k &       standard&    explicit\\
1 &       .9948&       .0899\\
2&        .0028   &    .0537\\
3 &       .0011&      .0502\\
4&        .0004 &      .0563\\
5 &       .0003&       .0696\\
6 &       .0002  &     .1860\\
7 &       .0001  &     .1107\\
8 &       .0001 &      .1194\\
9    &    .0001&       .1196\\
10 &     0. &   .1446
\end{tabular}
\end{center}
\end{table}
The histogram for the implicit filter is no longer perfectly balanced.
The asymmetry in the histogram reflects
the asymmetry of $F_0$ and can be eliminated by biasing $\xi$,
but there is no reason to do so; there is enough importance sampling without
this extra step. 

In Table V we display the estimates of the means of the density for the
two filters with 1000 particles for various values of $b$, compared with the exact results 
(the number of particles is relatively
large because with $h(x)=x^3$ and our parameter choices the variance of the conditional
density is significant, and this number of particles is needed for meaningful comparisons
of either algorithm with the exact result).   

\begin{table}
\begin{center}
Table V
\end{center}
Comparison of the the estimates of the means, implicit vs. standard filter,
$1000$ particles, together with the exact result, when $h(x)=x^3$,
as explained in the text. 
\begin{center}
\begin{tabular}{cccc}
b&        exact&       standard     &    implicit\\
0. &       0. &     -0.00 $\pm$.01 &    -.00 $\pm$.01\\
.5&         .109  &   .109 $\pm$.01 &    .109$\pm$.01\\
1.0   &      .442  &   .394 $\pm$ .04 &    .451$\pm$.02\\
1.5 &        .995  &   .775$\pm$.09  &    .995$\pm$.01\\
2.0  &      1.18   &   .875$\pm$.05  &   1.18$\pm$.01\\
2.5 &        1.30    &  .895 $\pm$.02 &   1.29$\pm$.02
\end{tabular}
\end{center}
\end{table}

As mentioned in the previous section, there are alternatives to the replacement of
$F$ by $F_0$; the point is that for each particle the function $F$ is an explicitly
known non-random function, and this fact can be used in multiple ways.

\section{Parameter identification}

One important application of particle filters is to parameter identification, where the SDE contains
an unknown parameter and the data are used to find this parameter's value. One of the standard ways
of doing this (see e.g \cite{Do2}) is system augmentation: 
one adds to the SDE the equation $d\sigma=0$  for the unknown parameter $\sigma$, 
one offers $\sigma$ a gamut of possible values, and one relies on the resampling process
that eliminate the values that do not fit the data. 
With the implicit filter this procedure fails, because the particles are not eliminated fast
enough. The alternative we are proposing is finding the unknown parameter $\sigma$ by stochastic approximation.
Specifically, 
Find a statistic $T$
of the output of the filter which is a function of $\sigma$, 
such that the 
expected value $E[T]$ vanishes when $\sigma$ has the right value $\sigma^*$,
and then solve the equation $E[T]=E[T(\sigma)]=0$ by the 
Robbins-Monro algorithm \cite{Rob}, in which the equation $E[T]=0$ is 
solved by the iteration:
\begin{equation}
\sigma_{n+1}=\sigma_n-\alpha_nT(\sigma_n),
\label{MM}
\end{equation}
where 
which converges when the coefficients $\alpha_n$ are such that $\sum \alpha_n \rightarrow
\infty$ while $\sum\alpha_n^2$ remains bounded.

As a concrete example, consider the SDE $dx=dW$,  
where $W$ is Brownian motion with variance $\sigma$, 
discretized with time steps $\delta$, with observations
$b^n=x^n+\eta$, where $\eta$ is a Gaussian with mean zero and variance
$s$. Data are generated by running the SDE once with the true  
value $\sigma^*$ of $\sigma$, adding the appropriate noise, and registering
the result at time $n\delta$ as $b^n$ for $n=1,2,\dots,N$.
For the functional $T$ we choose
\begin{equation}
T(\sigma)=C\sum (\Delta_i\Delta_{i-1})/\left((\sum \Delta_i^2)(\sum \Delta^2_{i-1})\right)^{1/2},
\end{equation}
where the summations are over $i$ between $2$ and $N$, $\Delta_i$
is the estimate of the increment of $x$ in the $i$-the step
and $C$ is a scaling constant. 
Clearly if the $\sigma$ used in the filtering equals $\sigma^*$
then by construction the successive values of $\Delta_i$ are independent
and $E[T]=0.$. 
We picked the parameters $N=100, \sigma=10^{-2}, s=10^{-4}, \delta=0.01$ (so that that
the increment of $W$ in one step has variance $10^{-4}$). 

Our algorithm is as follows: We make a guess $\sigma_1$, run the filter
for $N$ steps, evaluate $T$, and make a new guess for $\sigma$ using equation
(\ref{MM}) and $a_1=1$, rerun the filter, etc., with the $a_n$, the coefficient
in equation (\ref{MM}) at the $n$-th step, equal to $1/n$.
The scaling factor in (\ref{MM}) was found by trial and error: if it is too
large the iteration becomes unstable, if it is too small the
convergence is slow; we settled on $C=4$. 

This algorithm requires that the filter be run without either 
resampling or backward sampling, because resampling and backward
sampling introduce correlations between successive values of
the $\Delta_i$ and bias the values of $T$. In a long run, in particular in
a strongly nonlinear setting, one may need resampling for the filter to stay on track,
and this can be done by segmentation: divide the run of the filter into
segments of some moderate length $L$, perform the summations in the definition
of $T$ over that segment, then go back and run that segment with resampling,
then proceed to the next segment, etc.

The first question is,
how well is it possible in principle to reconstruct an unknown
value of $\sigma$ from $N$ observations; this issue was already
discussed in \cite{ct1}. Given $100$ samples of a Gaussian variable of
mean $0$ and variance $\sigma$, the variance reconstructed
from the observations is a random variable of mean $\sigma$ and variance $.16\cdot\sigma$;
$100$ observations do not contain enough information to reconstruct $\sigma$ perfectly.
A good way to estimate the best result
that can be achieved is to run the algorithm with the guess $\sigma_1$
equal to the exact value $\sigma^*$ with which the data were generated.
When this was done, the estimate of $\sigma$ was $1.27\sigma^*$.
This result indicates the achievable accuracy.

In Table VI we display the result of our algorithm when we start
with the starting value $\sigma_1=10\sigma^*$ and with $50$ particles. Each iteration
requires that one run the filter once.
\begin{table}
\begin{center}
Table VI
\end{center}
\centerline{Convergence of the parameter identification algorithm.}
\begin{center}
\begin{tabular}{cc}
Iteration   &  new estimate $\sigma/\sigma^*$ \\
\\
0  &            10.\\
1    &           .819\\
2      &         .943\\
3        &       1.02\\
4           &    1.05\\
5             &  1.08\\
6  &             1.10\\
7    &           1.13\\
8      &         1.15\\
9        &       1.16\\
10         &     1.17\\
11           &   1.18\\
12             &  1.18\\
13               & 1.18
\end{tabular}
\end{center}
\end{table}

\section{Conclusions}

We have presented the implicit filter for data assimilation,
together with several algorithms for the solution of the algebraic equations, 
including cases with non-convex functions $F$, as well as an algorithm for parameter
identification. The key idea in implicit sampling is to solve an algebraic
equation of the form $F(X)-\phi=\xi^T \xi/2$ for every particle, where the function $F$ is
explicitly known, $X$ is the new position of the particle, $\phi$ is an additive factor, and $\xi$
is a sample of a fixed reference pdf; $F$ varies from particle to particle and step to step. This construction makes it possible
to guide the particles to the high-probability area one by one under
a wide variety of circumstances. It is important to note that the equation
that links $\xi$ to $X$ is underdetermined and its solution can be adapted for
each particular problem.

Implicit sampling is of interest in particular because of its potential uses in high
dimensional problems, which are only briefly alluded to in the present paper.
The effectiveness of implicit sampling in high-dimensional settings depends on
one's ability to design maps $\xi \rightarrow x$ that satisfy the criteria above
and are computationally efficient. The design of such maps is problem dependent
and we will present examples in the context of specific applications. 

{\bf Acknowledgements} We would like to thank Prof. Jonathan Weare for asking penetrating
questions and for making very useful suggestions, Prof. Robert Miller for  
good advice and encouragement, and Mr. G. Zehavi for performing some of the preliminary computations. 
This work was supported in part by the Director, Office of Science,
Computational and Technology Research, U.S.\ Department of Energy under
Contract No.\ DE-AC02-05CH11231, and by the National Science Foundation under grants
DMS-0705910 and OCE-0934298.

\end{document}